\documentclass[a4paper,12pt]{article}

\newtheorem{theorem}{Theorem}
\newtheorem{lemma}{Lemma}

\newtheorem{remark}{Remark}

\usepackage{booktabs}

\usepackage{amsmath}
\usepackage{amssymb}
\usepackage{amsfonts}
\usepackage{mathrsfs} 
\usepackage{graphicx}
\usepackage{color}
\usepackage{ulem}

\newcommand{\new}[1]{{\color{blue}#1}}

\newcommand{\R}{\mathbb{R}} 
\newcommand{\X}{\mathbf{X}} 
\newcommand{\W}{\mathbf{W}} 
\renewcommand{\b}{\mathbf{b}} 
\newcommand{\A}{\mathbf{A}} 
\newcommand{\B}{\mathbf{B}} 

\title{\bf  Moment-SOS hierarchy and exit time of stochastic processes\footnote{
M. Velasco was partially supported by Proyecto INV-2018-50-1392 from Facultad de Ciencias, Universidad de los Andes. D. Henrion, M. Junca and M. Velasco were partially supported by ColCiencias Colombia-France cooperation Ecos Nord grant ``Problemas de momentos en control y optimizaci\'on'' (C19M02).
}}

\begin{document}

\author{Didier Henrion$^{1,2}$, Mauricio Junca$^3$, Mauricio Velasco$^3$}

\footnotetext[1]{CNRS, LAAS, Universit\'e de Toulouse, France. }
\footnotetext[2]{Faculty of Electrical Engineering, Czech Technical University in Prague, Czechia.}
\footnotetext[3]{Department of Mathematics, Universidad de los Andes, Bogot\'a, Colombia.}

\date{Draft of \today}

\maketitle

\begin{abstract}
The moment sum of squares (moment-SOS) hierarchy produces sequences of upper and lower bounds on functionals of the exit time solution of a polynomial stochastic differential equation with polynomial constraints, at the price of solving semidefinite optimization problems of increasing size. In this note we use standard results from elliptic partial differential equation analysis to prove convergence of the bounds produced by the hierarchy. We also use elementary convex analysis to describe a super- and sub-solution interpretation dual to a linear formulation on occupation measures. The practical relevance of the hierarchy is illustrated with numerical examples.
\end{abstract}

\section{Introduction}

This paper deals with the numerical evaluation of functionals of solutions of nonlinear stochastic differential equations (SDE).
Our approach consists of constructing a family of convex optimization problems (semidefinite programming problems, SDP) of increasing size whose solutions yield bounds on the given functional. This is an application of the so-called Lasserre or moment sum of squares (SOS) hierarchy \cite{l10,hkl20}. We are especially concerned about proving {\it convergence} of the bounds to the value of the functional.

The moment-SOS hierarchy was already developed and used in \cite{lp04,lpz06} for obtaining bounds on SDEs coming from finance. However, only lower and upper bounds were obtained, and the question of convergence was left open.

A key step to construct the moment-SOS hierarchy is the reformulation of the original, typically nonlinear problem, as a linear problem on occupation measures.
This linear reformulation is classical in Markov decision processes (MDP) \cite{bb96,ks98}. In order to prove convergence of the bounds obtained with the moment-SOS hierarchy, one has to prove that there is {\it no relaxation gap} between the original nonlinear problem and the linear problem on measures. This was already achieved in \cite{bb96,ks98} in the context of controlled MDP, but the proofs are lengthy and technical. Zero relaxation gap for optimal control of SDEs was proven in \cite{bgq11} with the help of viscosity solutions to Hamilton-Jacobi-Bellman partial differential equations (PDE). 

In \cite{kk13,k15}, bounds on functionals of solutions of SDEs were obtained by a dual approach, seeking test functions satisfying inequalities. When the functions and the SDE coefficients are polynomial, the inequalities are replaced by SOS constraints and solved numerically with SDP. Our occupation measure formulation can be interpreted as a primal approach, from which the dual on test functions follows from elementary convex analysis arguments. More recently, a primal-dual moment-SOS hierarchy approach to optimal control of SDEs was followed in \cite{jk20}, as a stochastic counterpart of \cite{hp17}, and no relaxation gap was ensured by approximating the value function solving the dual HJB PDE.

In this paper, we focus on a specific class of SDE functional evaluation, namely the {\it exit time of an uncontrolled SDE}. The exit time is a random variable that can be characterized by its moments. As shown in \cite{hrs01}, the exit time moments can be approximated numerically with occupation measures and linear programming (LP), with convergence guarantees based on the zero relaxation gap proof of \cite{ks98}. 

Our contribution is as follows:
\begin{itemize}
\item we provide a new proof of the equivalence, or zero relaxation gap, between the infinite-dimensional linear formulation on occupation measures and the original nonlinear SDE; the proof, much shorter and simpler in our opinion than the MDP proofs of \cite{bb96,ks98} or the HJB proof of \cite{bgq11,jk20}, relies on standard results from elliptic PDE analysis;
\item we describe a neat primal-dual linear formulation with no duality gap, allowing readily the application of the moment-SOS hierarchy.
\end{itemize}

The paper is organized as follows. The exit time problem is defined in Section \ref{sec:exit}. Its linear reformulation with occupation measures is described in Section \ref{sec:linear}. Our main result on zero relaxation gap is described and proved in Section \ref{sec:nogap}. The dual linear formulation is described in Section \ref{sec:dual}. Application of the moment-SOS hierarchy and numerical examples are described in Section \ref{sec:hierarchy}. Concluding remarks are gathered in Section \ref{sec:conclusion}

\section{Exit time problem}\label{sec:exit}

Let $\W(t)=(W_k)_{k=1,\ldots,m}$ denote the $m$-dimensional Brownian motion and let $\X(t)$ denote the solution of the stochastic differential equation (SDE)
\[
d\X = \b(\X)dt + \B(\X)d\W, \quad \X(0) = x
\] 
starting at $x \in {\mathscr X}$ where $\mathscr X$ is a given bounded open set of $\R^n$ with smooth boundary $\partial \mathscr X$ and closure $\overline{\mathscr X}:={\mathscr X} \cup {\partial\mathscr X}$. Drift functions $\b = (b_i)_{i=1,\ldots,n} : \R^n \to \R^n$ and diffusion functions $\B = (b_{ij})_{i=1,\ldots,n,\: j=1,\ldots,m} : \R^n \to \R^{n\times m}$ are given. We assume that $\B$ has full rank, so that the matrix $\A = (a_{ij}:=\frac{1}{2} \sum_{k=1}^m b_{ik}b_{jk})_{i,j=1,\ldots,n} : \R^n \to \R^{n\times n}$ is positive definite.
Assume $\b$ and $\B$ are continuous on $\R^n$ and growing at most linearly outside of $\mathscr X$, so that by standard arguments \cite[Chapter 5]{e13} there is a unique solution to the SDE, the stochastic process $\X(t)$.

Let $g : \partial \mathscr X \to \R$ be a given continuous function. We want to evaluate the function
\begin{equation}\label{value}
v^*(x) = E[g(\X(\tau_x))]
\end{equation}
where $\tau_x$ is the first time $\X(.)$ hits $\partial \mathscr X$, see e.g. \cite[Example 2, Section 6.2.1]{e13}.

\section{Linear reformulation}\label{sec:linear}

The generator of the stochastic process is the linear partial differential operator
\[
- L f : = \sum_{i,j=1}^n a_{ij} {\partial_i \partial_j f} + 
\sum_{i=1}^n b_i {\partial_i f}
\]
{where $\partial_i$ denotes the derivative with respect to the $i$-th variable}.
With this sign convention, and since the matrix $\A$ is positive definite, linear operator $L$ is uniformly elliptic \cite[Section 6.1.1]{e98}.

Given a function $f \in C^2(\overline{\mathscr X})$, It\^o's chain rule \cite[Chapter 4]{e13} implies that
\[
f(\X(\tau_x)) = f(\X(0)) - \int_0^{\tau_x} L f ds + \int_0^{\tau_x} Df \cdot \B \: d\W
\]
where
\[
Df \cdot \B \: d\W = \sum_{k=1}^m \sum_{i=1}^n {\partial_i f} b_{ik} \: dW_k.
\]
Taking the expected value yields Dynkin's formula  \cite[Section 6.1.3]{e13}:
\[
E[f(\X(\tau_x))] = \displaystyle E[f(\X(0))] - E\left[\int_0^{\tau_x} L f ds\right]
\]
that we rearrange as follows
\begin{equation}\label{dynkin}
E[f(\X(\tau_x))] + E\left[\int_0^{\tau_x} L f ds\right]  = \displaystyle f(x).
\end{equation}

Given $x \in \mathscr X$, define the expected occupation measure $\mu$ of the process $\X$ up to time $\tau_x$, such that 
\[
\mu({\mathscr A}) := E\left[\int_0^{\tau_x} I_{\mathscr A}(\X(s)) ds\right]
\]
for every set $\mathscr A$ in the Borel sigma algebra of $\mathscr X$, where $I_{\mathscr A}$ denotes the indicator function equal to one in $\mathscr A$ and zero outside. An equivalent analytic definition is
\[
{\langle f,\mu \rangle :=} E\left[\int_0^{\tau_x} f(\X(s))ds\right]
\]
for any test function $f$, {where
\[
\langle f,\mu \rangle := \int f \mu
\]
denotes the duality pairing of a continuous function $f$ and a
measure $\mu$}.
Define the exit location measure $\nu$ as the law of $\X({\tau_x})$ i.e.
\[
\nu({\partial\mathscr B}) := E[I_{\partial\mathscr B}(\X(\tau_x))]
\]
for every set $\mathscr B$ in the Borel sigma algebra of $\partial \mathscr X$. An equivalent analytic definition is
\[
{\langle f,\nu \rangle}
:= E\left[f(\X(\tau_x))\right]
\]
for every test function $f$.
Then Dynkin's formula \eqref{dynkin} becomes a linear partial differential equation on measures
\begin{equation}\label{linear}
{ \langle f,\nu \rangle + \langle L f,\mu \rangle} = f(x) 
\end{equation}
which can be written in the sense of distributions as
\begin{equation}\label{kolmo}
\nu + L' \mu = \delta_x
\end{equation}
where $L'$ is the linear operator adjoint to $L$
{and $\delta_x$ is the Dirac measure concentrated at $x$.} This equation is called the Kolmogorov or Fokker-Planck equation.

Following \cite{lp04,lpz06}, now define
\begin{equation}\label{minlp}
v_{\min}(x) := \min_{\mu,\nu} {\langle g,\nu \rangle} \:\:\mathrm{s.t.}\:\: \nu + L' \mu = \delta_x
\end{equation}
and
\begin{equation}\label{maxlp}
v_{\max}(x) := \max_{\mu,\nu} {\langle g, \nu \rangle} \:\:\mathrm{s.t.}\:\: \nu +   L' \mu = \delta_x
\end{equation}
which satisfy by construction
\[
v_{\min}(x) \leq v^*(x) \leq v_{\max}(x)
\]
for each $x \in \mathscr X$.
Note that  \eqref{minlp} and \eqref{maxlp} are linear optimization problems over measures $\mu$ and $\nu$ supported on $\mathscr X$ and $\partial \mathscr X$ respectively.

\section{No relaxation gap}\label{sec:nogap}

\begin{theorem}\label{nogap}
There is no relaxation gap between the nonlinear function evaluation problem \eqref{value} and the linear optimization problems   \eqref{minlp} and \eqref{maxlp}, i.e.
\[
v_{\min}(x)  = v^*(x) = v_{\max}(x)
\]
for each $x \in \mathscr X$.
\end{theorem}

{The proof of Theorem \ref{nogap} is based on the following result.}

{\begin{lemma}\label{unique}
For each $x \in \mathscr X$ there exists a unique exit location measure $\nu=\nu_x$ on $\partial\mathscr X$
solving the Kolmogorov equation \eqref{kolmo} for some expected occupation measure $\mu$ on $\mathscr X$ 
\end{lemma}}

{\bf Proof of Lemma \ref{unique}:}
Given any function $p \in C(\partial \mathscr X)$, let $f_p \in C^2(\mathscr X) \cup C(\partial \mathscr X)$ be the solution to the boundary value problem
\[
\begin{array}{cl}
L f = 0 &  \mathrm{in}\:\: \mathscr X \\
f = p & \mathrm{on}\:\:\partial \mathscr X
\end{array}
\]
which is unique according to \cite[Theorem 5, Section 2.2.3]{e98}. Plugging $f_p$ into \eqref{linear} yields 
\begin{equation}\label{nu}
{\new\langle p,\nu \rangle} = f_p(x).
\end{equation}
Since $\partial \mathscr X$ is compact, the space $C({\partial \mathscr X})$ is separable and by choosing countably many functions $p \in C({\partial \mathscr X})$ we can generate countably many linear relations \eqref{nu} that uniquely specify the measure $\nu$ that we denote $\nu_x$.
$\Box$

{
{\bf Proof of Theorem \ref{nogap}:}
Let $x \in \mathscr X$ and let $\nu_x$ denote the exit location measure on $\partial \mathscr X$ solving the Kolmogorov equation (\ref{kolmo}). Notice that the objective function in problems \eqref{minlp} and \eqref{maxlp} depends only on $\nu_x$. It follows that $\langle g, \nu_x \rangle = v_{\min}(x) = v_{\max}(x)$.
}
 
{
\begin{remark} It is natural to ask whether there is uniqueness of $\mu$ in Lemma \ref{unique}. The answer is clearly affirmative whenever the image of the operator $L$ is dense in the continuous functions on $\mathscr X$. For a specific instance when this occurs think of Brownian motion in the unit sphere in $\R^n$. The operator $L$ coincides with the Laplacian. The Laplacian maps homogeneous polynomials of degree $k$ surjectively onto homogeneous polynomials of degree $k-2$ and in particular it has dense image even when restricted to polynomials. Uniqueness of $\mu$ follows. By contrast if the image of $L$ is not dense then there are many non-trivial signed measures $\mu_0$ which annihilate $\langle Lf,\mu_0 \rangle$ for all $f$. For any such $\mu_0$ the pairs $(\nu_x, \mu+\mu_0)$ are solutions to the Kolmogorov equation \eqref{kolmo}. We do not know whether pairs of positive measures exist in all these cases.
\end{remark} 
} 
 
\section{Duality}\label{sec:dual}

{
In this section we use elementary notions from convex duality to derive the dual problem to the minimization resp. maximization problem on measures. We show that admissible solutions to the dual problem are subsolutions resp. supersolutions to the boundary value PDE solved by the value function. 
In particular we show that the concept of supersolution (resp. subsolution) arises naturally from elementary duality theory.
}

\begin{lemma}\label{dual}
The linear problem dual to \eqref{minlp} reads as follows
\begin{equation}\label{minlpdual}
\max_v v(x) \:\:\mathrm{s.t.}\:\: L v \leq 0 \:\:\mathrm{in}\:\:\mathscr X, \:\: v \leq g  \:\:\mathrm{on}\:\:\partial \mathscr X
\end{equation}
where the maximization is with respect to functions $v \in C^2({\mathscr X})$. There is no duality gap, i.e. the value of  \eqref{minlpdual} is equal to the value of \eqref{minlp}.
\end{lemma}

{\bf Proof:}
Let us denote by $v \in C(\mathscr X)$ the Lagrange multiplier corresponding to the equality constraint in primal problem \eqref{minlp}, and build the Lagrangian
$\ell(\mu,\nu,v) :=\langle g, \nu \rangle + \langle v,\delta_x-\nu-L'\mu \rangle = \langle g-v, \nu \rangle + \langle v,-L'\mu \rangle + \langle v, \delta_x \rangle = \langle g-v,\nu \rangle + \langle -Lv, \mu \rangle + v(x).$
The Lagrange dual function is then $\min_{\mu,\nu} \ell(\mu,\nu,v) = v(x)$ provided $v \leq g$ on $\partial\mathscr X$, the support of $\nu$, and $Lv \leq 0$ on $\mathscr X$, the support of $\mu$. The dual problem \eqref{minlpdual} then consists of maximizing the dual function subject to these inequality constraints. To prove that there is no duality gap, we use \cite[Theorem IV.7.2]{b02} and the fact that the image through the linear map $(\langle g,\nu\rangle, \: \nu + L'\mu)$ of the cone of measures $\mu$ resp. $\nu$ supported on $\mathscr X$ resp. $\partial\mathscr X$ is nonempty and bounded in the metric inducing the weak-star topology on measures.
$\Box$

As recalled in \cite[Example 2, Section 6.B]{e13}, the value function $v^*$ is the solution of the boundary value problem
\begin{equation}\label{bvp}
\begin{array}{cl}
L v = 0 &  \mathrm{in}\:\: \mathscr X \\
v = g & \mathrm{on}\:\:\partial \mathscr X.
\end{array}
\end{equation}

\begin{lemma}\label{subsol}
Any admissible function $v$ for linear problem \eqref{minlpdual} is a subsolution of boundary value problem \eqref{bvp}, in the sense that $v^* \geq v$ on $\overline{\mathscr X}$.
\end{lemma}

{\bf Proof:}
Let $v$ be admissible for \eqref{minlpdual}.
Function $u:=v-v^*$ is such that $Lu \leq 0$ in $\mathscr X$ and $u \leq 0$ on $\partial \mathscr X$.
By the weak maximum principle \cite[Theorem 1 page 327]{e98}, if $Lu \leq 0$ in $\mathscr X$ then $\max_{\overline{\mathscr X}} u = \max_{\partial \mathscr X} u$. Since $u \leq 0$ on $\partial \mathscr X$, this implies that $u \leq 0$ and hence $v^* \geq v$ on $\overline{\mathscr X}$.
$\Box$

Linear problem \eqref{minlpdual} selects the subsolution that touches the value function from below at $x$.

Similarly, the linear problem dual to \eqref{maxlp} reads as follows
\begin{equation}\label{maxlpdual}
\min_v v(x) \:\mathrm{s.t.}\: L v \geq 0 \:\:\mathrm{in}\:\:{\mathscr X}, \:\: v \geq g  \:\:\mathrm{on}\:\:\partial \mathscr X
\end{equation}
where the maximization is with respect to functions $v \in C^2({\mathscr X})$. There is no duality gap, i.e. the value of  \eqref{maxlpdual} is equal to the value of \eqref{maxlp}.
Any admissible function $v$ for linear problem \eqref{maxlpdual} is a super-solution of boundary value problem \eqref{bvp}, in the sense that $v^* \leq v$ on $X$ and $\partial \mathscr X$. Linear problem \eqref{maxlpdual} selects the supersolution that touches the value function from above at $x$.

\section{Random initial condition}

All the above developments generalize readily to the case that the initial condition $\X(0)$ in the SDE is a random variable whose law is a given probability measure $\xi$ on $\mathscr X$. The previous results can then be retrieved with the particular choice $\xi=\delta_x$ for a given $x \in \mathscr X$.

The quantity to be evaluated becomes
\[
v^*_\xi:=\int_{\mathscr X} E[g(\X(\tau_x))] d\xi(x).
\]
The Kolmogorov equation \eqref{kolmo}
becomes
\[
\nu + L' \mu = \xi
\]
and exactly the same arguments of the proof of Lemma \ref{unique} can be used to prove that it has a unique solution $\nu_\xi$ depending on $\xi$.
As in Theorem \ref{nogap}, it follows that the linear problems
\begin{equation}\label{primalmin}
\min_{\mu,\nu}\:\:\langle g,\nu \rangle \:\:\mathrm{s.t.}\:\: \nu + L' \mu = \xi
\end{equation}
and
\begin{equation}\label{primalmax}
\max_{\mu,\nu}\:\: \langle g,\nu \rangle \:\:\mathrm{s.t.}\:\: \nu + L' \mu = \xi
\end{equation}
have the same value $v^*_\xi$, i.e. there is no relaxation gap. The respective dual problems
\begin{equation}\label{dualmax}
\max_v \:\:\langle v,\xi \rangle \:\:\mathrm{s.t.}\:\: L v \leq 0 \:\:\mathrm{in}\:\:\mathscr X, \:\: v \leq g  \:\:\mathrm{on}\:\:\partial \mathscr X
\end{equation}
and
\begin{equation}\label{dualmin}
\min_v \:\:\langle v,\xi \rangle \:\:\mathrm{s.t.}\:\: L v \geq 0 \:\:\mathrm{in}\:\:\mathscr X, \:\: v \geq g  \:\:\mathrm{on}\:\:\partial \mathscr X
\end{equation}
have the same value $v^*_\xi$, i.e. there is no duality gap.

\section{Moment-SOS hierarchy and examples}\label{sec:hierarchy}

If the SDE coefficients $\b$ and $\B$ and the functional $g$ are semialgebraic\footnote{A function is semialgebraic if its graph is a semialgebraic set. A semialgebraic set is defined by a finite sequence of polynomial equations and inequalities.} in a semialgebraic set $\mathscr X$, we can apply the moment-SOS hierarchy on \eqref{minlp} and \eqref{maxlp} with convergence guarantees. In the primal, we obtain approximate moments (also called pseudo-moments or pseudo-expectations) of the occupation measures. They are not necessarily moments of the occupation measures as we are solving relaxations. In the SOS dual, we obtain polynomial sub- resp. super-solutions of increasing degrees of the boundary value problem \eqref{bvp}. {Each primal-dual problem is a semidefinite optimization problem.}

\subsection{Semidefinite relaxations}

 Let us briefly describe the construction of the moment-SOS hierarchy.
 
A bounded closed semialgebraic set ${\mathscr Z}$ of $\R^n$ can be written as the union of finitely many basic semialgebraic sets ${\mathscr Z}_i:=\{z\in \R^n \: : \: p_{i,j}(z) \geq 0, \:j=1,\ldots,m_i\}$, $i=1,\ldots,m$, described by finitely many polynomials $p_{i,j}$. Note that polynomial equations can be modeled by two inequalities of reverse signs. Since ${\mathscr Z}$ is bounded, without loss of generality, for each $i=1,\ldots,m$, one of the polynomials defining each ${\mathscr Z}_i$ can be chosen equal to $R^2-\sum_{k=1}^n z^2_k$ for $R$ sufficiently large, and for notational convenience we let $p_{i,0}(z):=1$
	
There are several algebraic characterizations of the set of positive polynomials on ${\mathscr Z}_i$. To describe one such characterization, let $r$ be a positive integer and for each $i=1,\ldots,m$ define the truncated quadratic module of degree $r$ of ${\mathscr Z}_i$, denoted $Q({\mathscr Z}_i)_r$, to be the set of polynomials which can be written as $\sum_{j=0}^{m_i} p_{i,j} s_{i,j}$ where the $s_{i,j}$ are sums of squares (SOS) of polynomials such that $2{\rm deg}s_{i,j}+{\rm deg}p_{i,j}\leq r$.
Every polynomial in the Minkowski sum $\sum_{i=1}^m Q({\mathscr Z}_i)_r$ is obviously nonnegative on ${\mathscr Z}$. Putinar's Positivstellensatz \cite{p93} is the much deeper statement that every polynomial strictly positive on ${\mathscr Z}$ lies in $\sum_{i=1}^m Q({\mathscr Z}_i)_r$.

We will now describe a hierarchy of semidefinite optimization problems which depend on the degree $r$ and provide us with upper and lower bounds on the value $v^*_\xi$. We will show that as $r\rightarrow \infty$ these bounds converge to the value $v^*_\xi$. 

Let us assume that $\mathscr X$ is a bounded basic semi-algebraic set. Its boundary $\partial \mathscr X$ is then a union of finitely many bounded basic semi-algebraic sets ${\mathscr X}^\partial_i$, $i=1,\ldots,m$. The primal (moment) problems are given by
\begin{equation}\label{mommin}
\begin{array}{rcl}
p^{\min}_r := & \min & \sum_{i=1}^{m} \ell_{\nu_i}(g) \\
& \mathrm{s.t.} & \ell_\mu(v) \geq 0, \quad \forall v \in Q({\mathscr X})_r \\
&&  \ell_{\nu_i}(v) \geq 0, \quad \forall v \in Q({\mathscr X}^\partial_i)_r \\
&& 
\ell_{\mu} (Lv)+\sum_{i=1}^m\ell_{\nu_i}(v)=\langle v,\zeta\rangle, \quad \forall v \in P_r
\end{array}
\end{equation}
and
\begin{equation}\label{mommax}
\begin{array}{rcl}
p^{\max}_r := & \max & \sum_{i=1}^{m} \ell_{\nu_i}(g) \\
& \mathrm{s.t.} & \ell_\mu(v) \geq 0, \quad \forall v \in Q({\mathscr X})_r \\
&&  \ell_{\nu_i}(v) \geq 0, \quad \forall v \in Q({\mathscr X}^\partial_i)_r \\
&& \ell_\mu(Lv)+\sum_{i=1}^m \ell_{\nu_i}(v) = \langle v,\xi \rangle, \quad \forall v \in P_r
\end{array}
\end{equation}
where the unknowns are linear operators $\ell_\mu, \ell_{\nu_1}, \ldots, \ell_{\nu_m}$ from $P_r$ to $\R$, for $P_r$ denoting the vector space of $n$-variate real polynomials of degree up to $r$. The dual (SOS) problems are given by
\begin{equation}\label{sosmax}
\begin{array}{rcl}
d^{\min}_r := & \max & \langle v,\xi \rangle \\
& \mathrm{s.t.} &  -Lv \in  Q({\mathscr X})_r \\
&& g-v \in \sum_{i=1}^m Q({\mathscr X}^\partial_i)_r
\end{array}
\end{equation}
and
\begin{equation}\label{sosmin}
\begin{array}{rcl}
d^{\max}_r := & \min & \langle v,\xi \rangle \\
& \mathrm{s.t.} &  Lv \in  Q({\mathscr X})_r \\
&& v-g \in \sum_{i=1}^m Q({\mathscr X}^\partial_i)_r
\end{array}
\end{equation}
where the unknowns are polynomials $v \in P_r$.

\begin{theorem}
Problems \eqref{mommin}, \eqref{mommax}, \eqref{sosmax} and \eqref{sosmin} are semidefinite programming problems. For each $r>0$, it holds
$d^{\min}_r \leq p^{\min}_r \leq v^*_\xi \leq p^{\max}_r \leq d^{\max}_r$. Moreover $\lim_{r\to\infty} (d^{\max}_r-d^{\min}_r)=0$.
\end{theorem}

{\bf Proof:}
To show that problems \eqref{sosmax} and \eqref{sosmin} are semidefinite programming problems, just observe that a polynomial $p(z)$ of degree at most $2r$ is a sum of squares of polynomials if and only if there is a positive semidefinite symmetric matrix $S$ such that $p(z)=b(z)'Sb(z)$ where $b(z)$ is a vector of polynomials spanning $P_r$. It follows that the truncated quadratic modules in \eqref{sosmax} and \eqref{sosmin} are projections of the semidefinite cone, and optimizing linear functions over them is an instance of semidefinite programming.

Problems \eqref{mommin} and \eqref{sosmax} are in duality. This follows easily from computing the Lagrangian as in Lemma~\ref{dual}. The constraint in problem \eqref{mommin} that a linear operator is non-negative for all polynomials in a truncated module can be expressed as a linear matrix inequality, i.e. it forms a spectrahedron, a linear slice of the semidefinite cone. From weak duality we conclude that $p^{\min}_r \geq d^{\min}_r$. Similarly, problems \eqref{mommax} and \eqref{sosmin} are in duality and $d^{\max}_r \geq p^{\max}_r$.

To prove that $p^{\min}_r \leq v^*_\xi$ observe that if $(\nu,\mu)$ are measures satisfying $\nu + \mathcal{L}'\mu=\xi$ then defining $\ell_\mu(v):=\langle v,\mu \rangle$ and $\ell_{\nu_i}(v):=\langle v,\nu_i \rangle$ where $\nu = \sum_{i=1}^{m} \nu_i$ for each $\nu_i$ supported on ${\mathscr X}^{\partial}_i$, we obtain admissible linear operators for the primal problem \eqref{mommin}.
Coefficients of the linear operators are moments of the respective measures.
It may however happen that linear operators admissible for problem \eqref{mommin} do not correspond to measures. Since we minimize a linear function on a possibly larger set, we obtain a lower bound. Similarly, we can prove that $v^*_\xi \leq p^{\max}_r$.

The most substantial claim is the convergence result. Given $\epsilon>0$ we will show that there exists an integer $r_\epsilon$ such that $0 \leq v^*_\xi-d^{\min}_r \leq \epsilon$ for $r\geq r_\epsilon$. Similar arguments imply an analogous convergence result for $d^{\max}_r$.

Recall from Section \ref{sec:dual} that the unique function $f\in C^2({\mathscr X})$ which satisfies $Lf=0$ in $\mathscr X$ and $f=g$ on $\partial \mathscr X$ is the value function $f=v^*$ of the problem and it satisfies $v^*_\xi=\langle v^*, \xi \rangle$. The proof proceeds by showing that this value function can be approximated by a sequence of elements $v_n$ which are feasible for problem \eqref{sosmax}. Since operator $L$ is uniformly elliptic, there is a polynomial $w$ such that $-Lw>0$ in $\mathscr X$, take e.g. a quadratic polynomial with sufficient large leading coefficients. By substracting a sufficiently large constant from $w$, we can also ensure that $w<g$ on $\partial \mathscr X$. Now let $\epsilon>0$ be given and let $(w_n)_{n\in {\mathbb N}}$ be a sequence of polynomials which approximate $v^*$ and its derivatives of order up to $2$ uniformly on $\overline{\mathscr X}$. Let $\eta>0$ be a real number with $\eta<\frac{\epsilon}{2\|w-v^*\|}$ and let $v_n:=(1-\eta)w_n+\eta w$. Due to the uniformity of the convergence of $w_n$ and the definition of $v_n$ the following statements hold for all $n$ sufficiently large: $\|w_n-v^*\|<\eta \|w-v^*\|$, $-Lv_n>0$ in $\overline{\mathscr X}$ and $v_n < g$ on $\partial \mathscr X$.
By Putinar's Positivstellensatz \cite{p93}, for any such $n$ there exists a degree $r_n$ such that 
$-Lv_n\in Q({\mathscr X})_{r_n}$ and $g-v_n\in Q({\partial \mathscr X})_{r_n}$ and therefore $v_n$ is feasible for problem \eqref{sosmax} for $r=r_n$. Furthermore, for any such $n$ we have
$\left|\langle v_n-v^*,\xi \rangle\right|\leq (1-\eta)\|w_n-v^*\|+\eta\|w-v^*\|\leq 2\eta\|w-v^*\|<\epsilon$
so we conclude that $v^*_\xi-d^{\min}_{r_n}\leq \epsilon$ proving the claim since $\epsilon>0$ was arbitrary.

\subsection{Solving the relaxations}

Moment relaxations of these linear problems can be modeled as generalized problems of moments with the GloptiPoly interface \cite{hll09} for Matlab. The relaxations are solved with the SDP solver in MOSEK \cite{d12}.

\subsection{Scalar example}\label{sec:ex1}

\begin{table}
\centering
\begin{tabular}{c|ccccccccc}
relaxation degree & 2 & 4 & 6 & 8 & 10 \\ \hline
lower bound & 0.65000 & 0.92157 & 0.98118 & 0.99503 & 0.99827 \\
upper bound & 1.00000 & 1.00000 & 1.00000 & 1.00000 & 1.00000  
\end{tabular}
\caption{Scalar example - bounds for increasing relaxation degrees.\label{tab:ex1}}
\end{table}

Let us illustrate the application of the moment-SOS hierarchy with an elementary exit time problem (\ref{value}) considered in \cite[Example 5.1]{hrs01}. The SDE is $d\X_t = (1+2\X_t)dt + \sqrt{2}\X_td\W_t$ on the domain $\mathscr X:=(0,1)$ with initial condition $x=1/2$. This process always exits at the point $\{1\}$, so $\nu=\delta_1$ solves \eqref{linear} and the value of linear problems \eqref{minlp} and \eqref{maxlp} is equal to $g(1)$. For the choice {$g(z)=z^2$} we report in Table \ref{tab:ex1} the values of the lower and upper bounds obtained by solving the moment relaxations to problems \eqref{minlp} and \eqref{maxlp}, for increasing relaxation degrees. The degree 10 relaxation was solved in 0.15 seconds on our laptop. The corresponding GloptiPoly script is given in the Appendix.

\subsection{Multivariate example}

\begin{figure}
  \centering
		\includegraphics[width=0.45\textwidth]{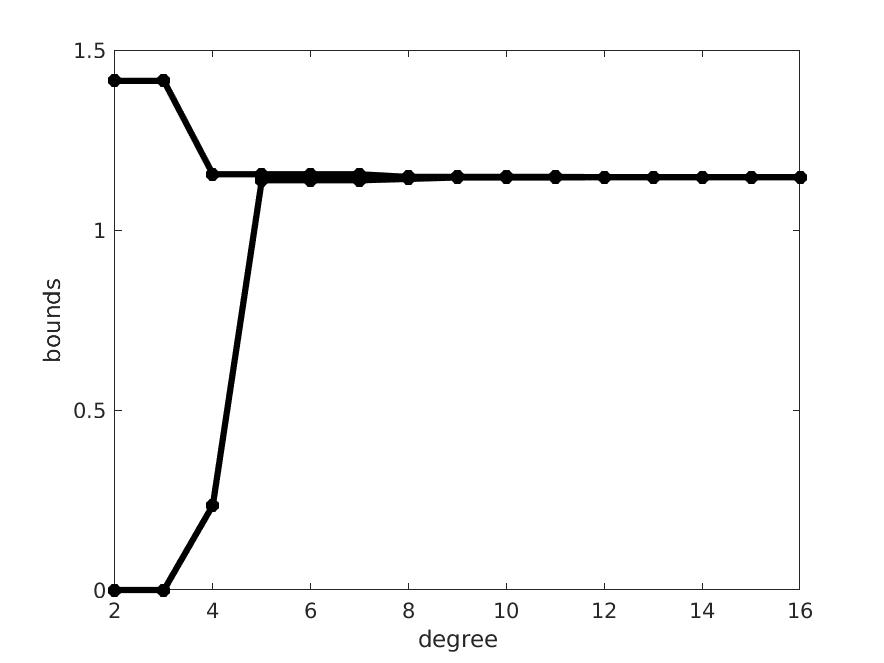}
	    \includegraphics[width=0.45\textwidth]{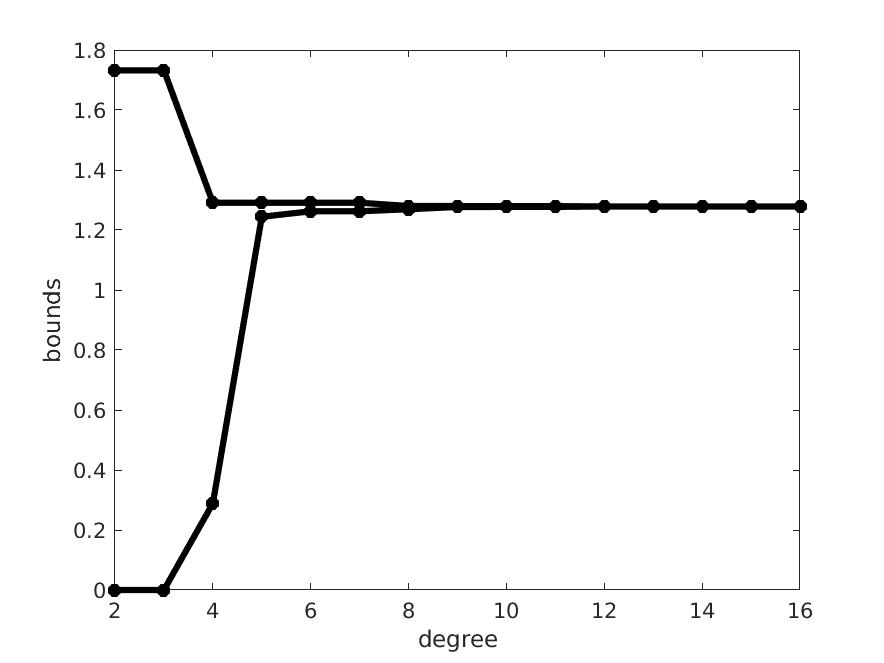}
	    \caption{Lower and upper bounds on the functional for increasing relaxation degrees and  dimension $n=2$ (left) and $n=3$ (right).\label{harmo}}
\end{figure}

Consider problem (\ref{value}) with {$g(z):=\sum_{k=1}^n z^2_k$} for the $n$-dimensional Brownian motion $\X_t=\W_t$, i.e. $\b=0$ and $\B=I_n$, in the convex semi-algebraic domain ${\mathscr X}:=\{{z \in \R^n : \sum_{k=1}^n z^4_k \leq 1}\}$ with initial condition $x=0$.

On Figure \ref{harmo} we plot for $n=2$ and $n=3$ the lower resp. upper bounds obtained by minimizing resp. maximizing the functional, for increasing relaxation degrees. We observe a fast convergence of the bounds.

In Table \ref{tab:ex2} we report the number of moments as well as the computational time required to solve the moment relaxation of degree $8$, for increasing values of the dimension $n$. For this relaxation degree the gap between the lower and upper bounds on the functional is less than 2\%. 

\begin{table}
	\centering
	\begin{tabular}{c|ccccccccc}
		dimension $n$ & 2 & 3 & 4 & 5 & 6 & 7 & 8 \\ \hline
		number of moments &  73 & 249 & 705 & 1749 & 3927 & 8151 & 15873\\
		CPU time (seconds) & 0.15 & 0.59 & 1.9 & 5.6 & 19 & 51 & 195\\
	\end{tabular}
	\caption{Computational burden for relaxation degree  8 and increasing dimensions.\label{tab:ex2}}
\end{table}

\section{Conclusion}\label{sec:conclusion}

Using elementary analytic arguments, we proved that there is no relaxation gap between the original problem and the linear problem on occupation measures in the special case of evaluating functionals of the exit time of stochastic processes on bounded domains. If the domain is basic semialgebraic and the SDE coefficients and the functional are polynomial or semi-algebraic, we can then readily apply the moment-SOS hierarchy with convergence guarantees. Tight bounds on the functionals can be obtained with off-the-shelf SDP solvers at a moderate cost.

{Of particular practical interest are approximations to the moments of the exit time distribution, as studied in \cite{hrs01}. In order to have access to these moments, the occupation measure and the boundary measure, as well as the test functions, should depend explicitly on time, as in \cite{lp04,lpz06}.}

We would also like to extend these techniques to optimal stopping time \cite{cs02} and stochastic optimal control problems \cite{bgq11}, with expectation constraints.

\section*{Appendix: Matlab script}\label{sec:script}

\begin{verbatim}
dmax = 10; % relaxation degree
mpol xmu xnu
mu = meas(xmu); % expected occupation measure
nu = meas(xnu); % exit location measure

x0 = 0.5; % initial condition
momeqs = []; % linear moment equations
for d = 0:dmax
 Lfmu = 0;
 if d > 0, Lfmu = Lfmu - mom((1+2*xmu)*(d*xmu^(d-1))); end
 if d > 1, Lfmu = Lfmu - mom(xmu^2*(d*(d-1)*xmu^(d-2))); end
 if d > 0, fnu = mom(xnu^d); else fnu = mass(nu); end
 momeqs = [momeqs; fnu+Lfmu == x0^d];
end

g = xnu^2; % functional

% construct moment relaxation
P = msdp(min(g), momeqs, xmu*(1-xmu)>=0, xnu*(1-xnu)==0);

% solve SDP problem
msol(P);

% bound
double(g)

% approximate mass of the occupation measure
double(mass(mu))
\end{verbatim}

\section*{Acknowledgement}

We are grateful to Milan Korda and Jean Bernard Lasserre for their feedback on this work.

\end{document}